\documentclass{article}

\usepackage[utf8]{inputenc} 		
\usepackage[english]{babel}		
\usepackage{amsmath,amssymb,amsthm}	
\usepackage{graphicx}			
\usepackage{hyperref}			
\usepackage{color}			

\newtheorem{theorem}{Theorem}[section]
\newtheorem{prop}[theorem]{Proposition}
\newtheorem{lemma}[theorem]{Lemma}

\theoremstyle{remark}

\newtheorem{rems}[theorem]{Remarks}

\title{Asymptotic study of supercritical surface Quasi-Geostrophic equation in critical space}
\author{Jamel BENAMEUR and Chaala KATAR\footnote{Department of Mathematics, Faculty of Science of Gab\`es, Research Laboratory Mathematics and Applications LR17ES11; Tunisia, jamelbenameur@gmail.com, katarchaala123@gmail.com}}

\begin{document}
\maketitle						

\tableofcontents					

\begin{center}
\textbf{Abstract} In this paper we prove, if $\theta\in C([0,\infty),H^{2-2\alpha}(\mathbb R^2))$ is a global solution of supercritical surface Quasi-Geostrophic equation with small initial data, then $\|\theta(t)\|_{H^{2-2\alpha}}$ decays to zero as time goes to infinity. Fourier analysis and
standard techniques are used.
\end{center}

\section{\bf Introduction}
The two-dimensional surface quasi-geostrophic equation with supercritical dissipation $(QG)$  is

   $$(QG)\hspace {2cm}\left\{\begin{array}{l}
   \partial_t \theta+|D|^{2\alpha}\theta +u_{\theta}.\nabla \theta =0\\
   \theta(0, x) =\theta^0(x).
   \end{array}\right.$$
where  $0< \alpha < 1/2$ a real number. The variable $\theta$ represents the potential temperature and  $u=(\partial_2|D|^{-1},\partial_1|D|^{-1} )\theta$ is the fluid velocity.
This equation is a two-dimensional model of the 3D incompressible Euler equations.
In the super critical case, there are only partial results.
We have global existence results when the initial data is small or local existence when the initial data are large.
The study of global solution is studied by serval researchers, See\cite{ChLee}, \cite{CC}, \cite{Ju}, \cite{Miu}, ...\\
The equation is invariant under the following scaling:
 $$\theta_\lambda(t,x)=\lambda^{2\alpha-1}\theta(\lambda^{2\alpha}t,\lambda x)$$ with initial data $\theta^0_{\lambda}= \lambda^{2\alpha-1} \theta^0(\lambda x)$, so $\dot{H}^{2-2\alpha}$ is a critical space and we have
 $$\|\lambda^{2\alpha-1} f(\lambda.)\|_{\dot{H}^{2-2\alpha}}=\|f\|_{\dot{H}^{2-2\alpha}}, \forall \lambda>0$$

 We begin by recalling the uniqueness and existence result in critical space.
\begin{theorem}\label{th00}(\cite{HD})
 Let $\theta^0 \in H^{2-2\alpha}(\mathbb{R}^2)$ . Then, there is a unique time $T^*\in(0, \infty]$ and a unique solution $\theta \in C([0, T^*), H^{2-2\alpha}(\mathbb{R}^2))$ of $(QG)$.
\end{theorem}
Our work is devoted to study the maximal solution of $(QG)$ system given by Theorem \ref{th00}. First, we study the case of a non-regular solution, where we show that we have a blow up at the first time of irregularity. Secondly, the case of a regular solution, we show that if moreover the initial condition is small enough then the solution decreases towards 0 at infinity. The proof of explosion result is classic and similar to that of the incompressible Navier-Stokes equations, the only thing to point out that one needs to specify some estimate and the dependence of the associated constants of the index of parametric regularity ($C=\sigma2^\sigma C(\alpha)$). But the novelty are the tools used to show the result on the regular solutions where one has to go down in the index of regularity of type $\dot H^{-\sigma}$ for the high frequencies. Our first main result reads as follows:
\begin{theorem}\label{TEC1}
Let $\theta\in C([0,T^*),H^{2-2\alpha}(\mathbb{R}^2))$ be a maximal solution of $(QG)$. If $T^*<\infty$, then
$$\int^{T^*}_0 \||D|^{\alpha}\theta(\tau)\|^2_{H^{2-2\alpha}}d\tau =\infty.$$
\end{theorem}
\begin{theorem}\label{TEC2} There is a positif real number $\varepsilon_0$ such that: If $\theta^0 \in H^{2-2\alpha}(\mathbb{R}^2)$ verifies $\|\theta^0\|_{H^{2-2\alpha}}<\varepsilon_0$, then  there exists a global solution of $(QG)$ such that $$\theta \in C([0,\infty),H^{2-2\alpha}(\mathbb{R}^2))\cap L^2(\mathbb R^+,\dot H^{2-\alpha}(\mathbb{R}^2)).$$
In addition the following estimate holds
\begin{equation}\label{eq1th1}\|\theta\|_{H^{2-2\alpha}}^2+\int_0^t\||D|^\alpha\theta\|_{H^{2-2\alpha}}^2\leq \|\theta^0\|_{H^{2-2\alpha}}^2.\end{equation}
\end{theorem}

\begin{theorem}\label{TEC3} If $\theta \in C([0,\infty),H^{2-2\alpha}(\mathbb{R}^2))$ is a global solution of $(QG)$ such that $ \|\theta(0)\|_{H^{2-2\alpha}}<\varepsilon_0$, then
$$\lim_{t\rightarrow\infty}\|\theta(t)\|_{H^{2-2\alpha}}=0.$$
\end{theorem}
\begin{rems} The proof of Theorem \ref{TEC3} is done in two steps. In the first, we prove $\|\theta(t)\|_{L^2}\rightarrow0$ if $t\rightarrow\infty$, where we use the $L^2$-energy estimate
$$\|\theta(t)\|_{L^2}^2+2\int_0^t\||D|^\alpha \theta(z)\|_{L^2}dz\leq \|\theta^0\|_{L^2}^2$$
and the uniqueness given by Theorem \ref{th00}. In the second step, we prove $\|\theta(t)\|_{\dot H^{2-2\alpha}}\rightarrow0$ if $t\rightarrow\infty$, where we use
the uniqueness in $H^{2-2\alpha}(\mathbb R^2)$ and $$\theta\in L^\infty([0,\infty),L^2(\mathbb R^2))\cap L^2([0,\infty),\dot H^{2-\alpha}(\mathbb R^2)).$$
\end{rems}

\noindent The rest of this paper is organized as follows: In Section 2, we recall some notations and we give some important preliminary results. In section 3, we prove a blow-up criterion to the non regular solution of $(QG)$ in critical space. Section 4 is devoted to give a proof of Theorem \ref{TEC1}. In section 5, we prove that the norm of global solution in $H^{2-2\alpha}$ goes to zero at infinity, which is made in two steps.
\section{\bf Notations and preliminaries results}
\subsection{Notations}
In this short section, we collect some notations and definitions that will be used later.
\begin{enumerate}
\item[$\bullet$] If $A$ is a subset of $\mathbb R^2$, ${\bf 1}_{A}$ denotes the characteristic function of $A$.
\item[$\bullet$] The Fourier transformation is normalized as
$$\mathcal{F}(f)(\xi)=\widehat{f}(\xi)=\int_{\mathbb R^2}\exp(-ix.\xi)f(x)dx,\,\,\,\xi=(\xi_1,\xi_2)\in\mathbb R^2.$$
\item[$\bullet$] The inverse Fourier formula is
$$\mathcal{F}^{-1}(g)(x)=(2\pi)^{-2}\int_{\mathbb R^3}\exp(i\xi.x)g(\xi)d\xi,\,\,\,x=(x_1,x_2)\in\mathbb R^2.$$
\item[$\bullet$] The convolution product of a suitable pair of function $f$ and $g$ on $\mathbb R^2$ is given by
$$(f\ast g)(x):=\int_{\mathbb R^2}f(y)g(x-y)dy.$$
\item[$\bullet$] If $f:\mathbb R^2\rightarrow\mathbb R$ is a function and $v=(v_1,v_2,v_3)$ is a vector fields, we set
$${\rm div}\,(fv):=v.(\nabla f)+f{\rm div}\,(v).$$
Moreover, if $\rm{div}\,v=0$ we obtain
$${\rm div}\,(fv):=v.(\nabla f).$$
\item[$\bullet$] If $f\in\mathcal S'(\mathbb R^2)$ such that $\widehat{f}\in L^1_{loc}(\mathbb R^2)$, we define $|D|^\sigma f$ by
$$\widehat{|D|^\alpha f}(\xi)= |\xi|^\alpha\widehat{f}(\xi).$$
\item[$\bullet$] For $s\in \mathbb{R}, H^s (\mathbb{R}^2)$ denotes the usual non-homogeneous Sobolev space on $\mathbb{R}^2$
and $\langle.,.\rangle_{H^s}$ denotes the usual scalar product on $H^s(\mathbb{R}^2)$.
\item[$\bullet$] For $s\in \mathbb{R}, \dot{H}^s (\mathbb{R}^2)$ denotes the usual homogeneous Sobolev space on $\mathbb{R}^2$
and $\langle.,.\rangle_{\dot{H}^s}$ denotes the usual scalar product on $\dot{H}^s(\mathbb{R}^2)$.
\item[$\bullet$] For $s>0$, we have $(1+|\xi|^2)^s\sim 1+|\xi|^{2s}$, then $$\|f\|_{H^s}^2\sim\|f\|_{L^2}^2+\||D|^sf\|_{L^2}^2.$$
We can therefore use the equivalent norm $$\sqrt{\|f\|_{L^2}^2+\||D|^sf\|_{L^2}^2}$$ and the associated scalar product. This scalar product helps us use the property $$\langle v.\nabla f/f\rangle_{L^2}=0\;\;{\rm if}\;\;{\rm div}\;(v)=0,$$ and$$\langle v.\nabla f/f\rangle_{H^s}=\langle v.\nabla f/f\rangle_{\dot H^s}.$$
\end{enumerate}
\subsection{Preliminaries results}
In this section we recall some classical results and we present new technical lemmas.
\begin{prop}(\cite{HBAF})\label{thhb} Let $H$ be Hilbert space and $(x_n)$ be a bounded sequence of elements in $H$ such that
$$x_n\rightarrow x\;weakly\;in\;H\;\;and\;\;\limsup_{n\rightarrow\infty}\|x_n\|\leq \|x\|,$$
then $\lim_{n\rightarrow\infty}\|x_n-x\|=0.$
\end{prop}
\begin{lemma}(\cite{JYC})\label{LP}
Let $s_1,\ s_2$ be two real numbers such that $s_1<1\; and\; s_1+s_2>0$. There exists a constant  $C_1=C_1(s_1,s_2),$ , such that for all $f,g\in \dot{H}^{s_1}(\mathbb{R}^2)\cap \dot{H}^{s_2}(\mathbb{R}^2),$ $f.g \in \dot{H}^{s_1+s_2-1}(\mathbb{R}^2)$ and
 \begin{equation}\label{IS1}
          \|fg\|_{\dot{H}^{s_1+s_2-1}}\leq C_1 (\|f\|_{\dot{H}^{s_1}}\|g\|_{\dot{H}^{s_2}}+\|f\|_{\dot{H}^{s_2}}\|g\|_{\dot{H}^{s_1}}).
 \end{equation}
          If in addition $s_2<1,$ there exists a constant $C_2=C_2(s_1,s_2)$ such that for all $f \in \dot{H}^{s_1}(\mathbb{R}^2)$ and $g\in\dot{H}^{s_2}(\mathbb{R}^2)$, then  $f.g \in \dot{H}^{s_1+s_2-1}(\mathbb{R}^2)$ and
          \begin{equation}\label{IS2}
          \|fg\|_{\dot{H}^{s_1+s_2-1}}\leq C_2 \|f\|_{\dot{H}^{s_1}}\|g\|_{\dot{H}^{s_2}}.
 \end{equation}
 \end{lemma}
\begin{rems}
Let $\alpha\in(0,1/2)$, then there is a constant $C(\alpha)$ such that:\\
If $f\in \dot H^0(\mathbb R^2)=L^2(\mathbb R^2)$ and $g\in \dot H^{1-\alpha}(\mathbb R^2)$, then $fg\in\dot H^{-\alpha}(\mathbb R^2)$ and
$$\|fg\|_{\dot H^{-\alpha}}\leq C(\alpha)\|f\|_{L^2}\|g\|_{\dot H^{1-\alpha}}.$$
\end{rems}
 \begin{lemma}\label{01}
          Let $0<\alpha<1/2$, there is a constant $C(\alpha)$ such that for $\theta \in H^{\sigma+\alpha}(\mathbb{R}^2)\cap H^{2-2\alpha}(\mathbb{R}^2)$ with $\sigma \geq 1$, we have
           $$|\langle u_{\theta}.\nabla\theta, \theta \rangle_{H^{\sigma}}|\leq \sigma2^\sigma C(\alpha)\|\theta\|_{\dot{H}^{2-2\alpha}} \|\theta\|^2_{\dot{H}^{\sigma+\alpha}}.$$
 \end{lemma}
{\bf Proof.} Using the fact $|\langle u_{\theta}.\nabla\theta, \theta \rangle_{L^2}=0$, (${\rm div\,} u_{\theta}=0$),  we get
$$ |\langle u_{\theta}.\nabla\theta, \theta \rangle_{H^\sigma}|\leq \int_{\xi}\int_{\eta}||\xi|^\sigma-|\eta|^\sigma||\widehat{u_\theta}(\xi-\eta)||\widehat{\nabla\theta}(\eta)|d\eta|\xi|^\sigma|\widehat{\theta}(\xi)|d\xi.$$
By using the elementary inequality
$$||\xi|^\sigma-|\eta|^\sigma|\leq \sigma 2^{\sigma-1}|\xi-\eta|(|\eta|^{\sigma-1}+|\xi-\eta|^{\sigma-1})$$
we get          $$\begin{array}{ll}
          |\langle u_{\theta}.\nabla\theta, \theta \rangle_{H^\sigma}|&\leq \displaystyle|\int_{\xi}\int_{\eta}|(|\xi|^\sigma-|\eta|^\sigma)|.|\widehat{\theta}(\xi-\eta)||\widehat{\nabla\theta}(\eta)|d\eta|\xi|^\sigma|\widehat{\theta}(\xi)|d\xi|\\
         &\leq \displaystyle \sigma2^{\sigma-1}\int_{\xi}\int_{\eta}|\xi-\eta|^\sigma|\widehat{\theta}(\xi-\eta)|.|\eta||\widehat{\theta}(\eta)|d\eta|\xi|^\sigma|\widehat{\theta}(\xi)|d\xi\\
         &\displaystyle+\sigma2^{\sigma-1}\int_{\xi}\int_{\eta}|\xi-\eta||\widehat{\theta}(\xi-\eta)|.|\eta|^\sigma|\widehat{\theta}(\eta)|d\eta|\xi|^\sigma|\widehat{\theta}(\xi)|d\xi\\
          &\leq \displaystyle\sigma2^{\sigma-1}\Big( \int_{\xi}|\xi|^{-2\alpha}(\int_{\eta}|\xi-\eta|^\sigma|\widehat{\theta}(\xi-\eta)||\widehat{\nabla\theta}(\eta)|d\eta)^2d\xi\Big)^{1/2} \|\theta\|_{\dot{H}^{s+\alpha}}\\
          &+\displaystyle\sigma2^{\sigma-1}\Big(\int_{\xi}|\xi|^{-2\alpha}(\int_{\eta}|\xi-\eta||\widehat{\theta}(\xi-\eta)||\eta|^\sigma|\widehat{\theta}(\eta)|d\eta)^2d\xi\Big)^{1/2} \|\theta\|_{\dot{H}^{\sigma+\alpha}}\\
          &\leq \sigma2^{\sigma}\|fg\|_{\dot H^{-\alpha}}\|\theta\|_{\dot{H}^{\sigma+\alpha}},
          \end{array}$$
with
$$\left\{\begin{array}{lcl}
\widehat{f}&=&|\xi|.|\widehat{\theta}(\xi)|\\
\widehat{g}&=&|\xi|^\sigma|\widehat{\theta}(\xi)|.
\end{array}\right.$$
Using the product laws, with $ s_1+s_2=1-\alpha>0$ and
$$\left\{\begin{array}{l}
s_1=1-2\alpha<1\\
s_2=\alpha<1
\end{array}\right.$$
 we obtain the desired result.
\begin{lemma}\label{PS}
Let $0<\alpha<1/2$ and for $\theta,\omega \in H^{2-\alpha}(\mathbb{R}^2)$ we have
\begin{equation}\label{eq23.1}|\langle u_{\omega}.\nabla\theta, \theta \rangle_{H^{2-2\alpha}}|\leq C(\alpha)\|\omega\|_{\dot{H}^{2-\alpha}} \|\theta\|_{\dot{H}^{2-2\alpha}}\|\theta\|_{\dot{H}^{2-\alpha}},\end{equation}
\begin{equation}\label{eq23.2}|\langle u_{\omega}.\nabla\theta, \theta \rangle_{H^{2-2\alpha}}|\leq C(\alpha)\|\omega\|_{\dot{H}^{2-2\alpha}}\|\theta\|_{\dot{H}^{2-\alpha}}^2.\end{equation}
\end{lemma}
{\bf Proof.} Proof of equation (\ref{eq23.1}): Put $\delta=2-2\alpha$, we have $$\begin{array}{ll}
|\langle u_{\omega}.\nabla\theta, \theta \rangle_{H^{\delta}}|&=|\langle u_{\omega}.\nabla\theta, \theta \rangle_{\dot H^{\delta}}|\\
&=|\langle|D|^\delta (u_{\omega}.\nabla\theta), |D|^\delta\theta \rangle_{L^2}|\\
&=|\langle|D|^\delta (u_{\omega}.\nabla\theta), |D|^\delta\theta \rangle_{L^2}-\langle u_{\omega}\nabla.|D|^\delta\theta, |D|^\delta\theta \rangle_{L^2}|\\
&\leq\displaystyle C\Big( \int_{\xi}|\xi|^{-2\alpha}(\int_{\eta}|\xi-\eta|^\delta|\widehat{\omega}(\xi-\eta)||\widehat{\nabla\theta}(\eta)|d\eta)^2d\xi\Big)^{1/2} \|\theta\|_{\dot{H}^{\delta+\alpha}}\\
          &\displaystyle +C\Big(\int_{\xi}|\xi|^{-2\alpha}(\int_{\eta}|\xi-\eta||\widehat{\omega}(\xi-\eta)||\eta|^\delta|\widehat{\theta}(\eta)|d\eta)^2d\xi\Big)^{1/2} \|\theta\|_{\dot{H}^{\delta+\alpha}}\\
          &\leq C(\|f_2g_2\|_{\dot H^{-\alpha}}+\|f_3g_3\|_{\dot H^{-\alpha}})\|\theta\|_{\dot{H}^{\delta+\alpha}}.
\end{array}$$
with
$$\left\{\begin{array}{lcl}
\widehat{f_2}&=&|\xi|^\delta.|\widehat{\omega}(\xi)|\\
\widehat{g_2}&=&|\xi|.|\widehat{\theta}(\xi)|\\
\widehat{f_3}&=&|\xi|.|\widehat{\omega}(\xi)|\\
\widehat{g_3}&=&|\xi|^\delta|\widehat{\theta}(\xi)|.
\end{array}\right.$$
Using the product laws, with
$${\rm for}\;(f_2,g_2):\;\left\{\begin{array}{l}
s_1+s_2=1-\alpha>0\\
s_1=\alpha<1\\
s_2= 1-2\alpha <1
\end{array}\right.\;\;\;$$
and
$${\rm for}\;(f_3,g_3):\;\left\{\begin{array}{l}
s_1+s_2=1-\alpha>0\\
s_1=1-\alpha <1\\
s_2=0<1.
\end{array}\right.$$
we obtain the desired result.\\
Proof of equation (\ref{eq23.2}): We have $$\begin{array}{ll}
|\langle u_{\omega}.\nabla\theta, \theta \rangle_{H^{\delta}}|&=|\langle u_{\omega}.\nabla\theta, \theta \rangle_{\dot H^{\delta}}|\\
&=|\langle|D|^\delta (u_{\omega}.\nabla\theta), |D|^\delta\theta \rangle_{L^2}|\\
&=|\langle|D|^\delta (u_{\omega}.\nabla\theta), |D|^\delta\theta \rangle_{L^2}-\langle u_{\omega}\nabla.|D|^\delta\theta, |D|^\delta\theta \rangle_{L^2}|\\
&\leq\displaystyle  C\Big( \int_{\xi}|\xi|^{-2\alpha}(\int_{\eta}|\xi-\eta|^\delta|\widehat{\omega}(\xi-\eta)||\widehat{\nabla\theta}(\eta)|d\eta)^2d\xi\Big)^{1/2} \|\theta\|_{\dot{H}^{\delta+\alpha}}\\
          &\displaystyle +C\Big(\int_{\xi}|\xi|^{-2\alpha}(\int_{\eta}|\xi-\eta||\widehat{\omega}(\xi-\eta)||\eta|^\delta|\widehat{\theta}(\eta)|d\eta)^2d\xi\Big)^{1/2} \|\theta\|_{\dot{H}^{\delta+\alpha}}\\
          &\leq C(\|f_4g_4\|_{\dot H^{-\alpha}}+\|f_5g_5\|_{\dot H^{-\alpha}})\|\theta\|_{\dot{H}^{\delta+\alpha}}.
\end{array}$$
with
$$\left\{\begin{array}{lcl}
\widehat{f_4}&=&|\xi|^\delta.|\widehat{\omega}(\xi)|\\
\widehat{g_4}&=&|\xi|.|\widehat{\theta}(\xi)|\\
\widehat{f_5}&=&|\xi|.|\widehat{\omega}(\xi)|\\
\widehat{g_5}&=&|\xi|^\delta|\widehat{\theta}(\xi)|.
\end{array}\right.$$
Using the product laws, with
$${\rm for}\;(f_4,g_4):\;\left\{\begin{array}{l}
s_1+s_2=1-\alpha>0\\
s_1=0<1\\
s_2= 1-\alpha <1
\end{array}\right.$$
and
$${\rm for}\;(f_5,g_5):\;\left\{\begin{array}{l}
s_1+s_2=1-\alpha>0\\
s_1=1-2\alpha <1\\
s_2=\alpha<1,
\end{array}\right.$$
we obtain the desired result.
\begin{lemma}\label{lemf} Let $h:[0,T]\rightarrow[0,\infty]$ be measurable function and $\sigma>0$, then
$$\Big(\int_0^Te^{-\sigma(T-z)}h(z)dz\Big)^2\leq 2\sigma^{-1}\int_0^Te^{-\sigma(T-z)}h(z)^2dz.$$
\end{lemma}
{\bf Proof.} We begin by denote $I=\displaystyle\Big(\int_0^Te^{-\sigma(T-z)}h(z)dz\Big)^2$. We have
$$I=\displaystyle\int_0^T\int_0^T e^{-\sigma(T-z)}e^{-\sigma(T-z')}h(z)h(z')dzdz'.$$
Using the fact $h(z)h(z')\leq \sup(h(z)^2,h(z')^2)\leq h(z)^2+h(z')^2$, we get
$$\begin{array}{lcl}
I&\leq&\displaystyle\int_0^T\int_0^T e^{-\sigma(T-z)}h(z)^2e^{-\sigma(T-z')}dzdz'+\int_0^T\int_0^T e^{-\sigma(T-z)}h(z')^2e^{-\sigma(T-z')}dzdz'\\
&\leq&\displaystyle(\int_0^Te^{-\sigma(T-z')}dz').\int_0^T e^{-\sigma(T-z)}h(z)^2dz+(\int_0^Te^{-\sigma(T-z)}dz)\int_0^T h(z')^2e^{-\sigma(T-z')}dz'\\
&\leq&\displaystyle(\frac{1-e^{-\sigma T}}{\sigma})\int_0^T e^{-\sigma(T-z)}h(z)^2dz+(\frac{1-e^{-\sigma T}}{\sigma})\int_0^T e^{-\sigma(T-z)}h(z)^2dz\\
&\leq&\displaystyle\frac{2}{\sigma}\int_0^T e^{-\sigma(T-z)}h(z)^2dz.
\end{array}$$
\section{\bf Blow up of non regular solution}
In this section we prove Theorem \ref{TEC1}. This proof is done in three steps. Let $\theta \in C([0,T^*),H^{2-2\alpha}(\mathbb{R}^2))$ be a maximal solution of the system $(QG)$. We suppose that $T^*<\infty$ and \begin{equation}\label{bleq0}\int^{T^*}_0 \| |D|^{\alpha}\theta\|^2_{H^{2-2\alpha}}d\tau <+\infty.\end{equation}
{\bf Step 1:} We prove that $\theta$ is bounded in $H^{2-2\alpha}(\mathbb R^2)$. For this select a time $T_0 \in (0, T^*)$ such that $\displaystyle\int^{T^*}_{T_0} \| |D|^{\alpha}\theta\|^2_{H^{2-2\alpha}}d\tau <\frac{1}{4C(\alpha)}.$ Then, by using (\ref{eq23.1}) we get: for all $t\in [T_0,T^*)$
$$\begin{array}{ll}
\displaystyle \|\theta(t)\|^2_{{H}^{2-2\alpha}}+ 2\int^t_{T_0}\|\theta\|^2_{{H}^{2-\alpha}}d\tau&\displaystyle \leq \|\theta(T_0)\|^2_{{H}^{2-2\alpha}}+2\int^t_{T_0}|<u_\theta(\tau).\nabla \theta(\tau)/\theta(\tau)>_{\dot{H}^{2-2\alpha}}|d\tau\\
&\leq \displaystyle \|\theta(T_0)\|^2_{{H}^{2-2\alpha}}+2C(\alpha)\int^t_{T_0}\|\theta(\tau)\|_{\dot{H}^{2-2\alpha}}\|\theta\|^2_{\dot{H}^{2-\alpha}}d\tau\\
&\displaystyle\leq\|\theta(T_0)\|^2_{{H}^{2-2\alpha}}+2C(\alpha)\sup_{z\in[T_0,t]}\|\theta(z)\|_{\dot{H}^{2-2\alpha}}\int^t_{T_0}\|\theta\|^2_{\dot{H}^{2-\alpha}}d\tau\\
&\displaystyle \leq\|\theta(T_0)\|^2_{{H}^{2-2\alpha}}+\frac{1}{2}\sup_{z\in[T_0,t]}\|\theta(z)\|_{\dot{H}^{2-2\alpha}}.
\end{array}$$
Then $$\|\theta(t)\|^2_{{H}^{2-2\alpha}}\leq 2\|\theta(T_0)\|^2_{{H}^{2-2\alpha}}, \quad \forall t\in [T_0,T^*)$$
Put $$ M=\max(\sup_{z\in[0,T_0]}\|\theta(z)\|_{{H}^{2-2\alpha}}; \sqrt{2}\|\theta(T_0)\|_{{H}^{2-2\alpha}}).$$
{\bf Step 2:} We want to prove $\theta$ is a Cauchy family at $T^*$ in $L^2(\mathbb{R}^2)$. For $t<t'\in[T_0,T^*)$, we have
$$\theta(t)-\theta(t')=-\int^{t'}_t|D|^{2\alpha}\theta-\int^{t'}_t u_\theta\nabla\theta.$$
Then
$$\|\theta(t)-\theta(t')\|_{L^2}\leq\int^{t'}_t\||D|^{2\alpha}\theta-u_\theta\nabla\theta\|_{L^2}.$$
By using Lemma \ref{LP} with $s_1=1-2\alpha$ and $s_2=2\alpha$, we get
$$\begin{array}{ll}
\displaystyle \|\theta(t)-\theta(t')\|_{L^2}&\displaystyle\leq\int^{t'}_t\|\theta\|_{H^{2\alpha}}+\int^{t'}_t \|u_\theta\nabla\theta\|_{L^2}\\
&\leq \displaystyle \int^{t'}_t\|\theta\|_{\dot{H}^{2\alpha}}+C(\alpha)\int^{t'}_t\|\theta\|_{\dot{H}^{2\alpha}}\|\theta\|_{\dot{H}^{2-2\alpha}}\\
&\leq \displaystyle \int^{t'}_t\|\theta\|_{\dot H^{2\alpha}}+C(\alpha)M\int^{t'}_t\|\theta\|_{\dot{H}^{2\alpha}}\\
&\leq \displaystyle (1+C(\alpha)M) \int^{t'}_t\|\theta\|_{\dot{H}^{2\alpha}}\\
&\leq \displaystyle (1+C(\alpha)M) \int^{t'}_t\|\theta\|_{{H}^{2\alpha}}\\
&\leq \displaystyle (1+C(\alpha)M)\int^{t'}_t\|\theta\|_{{H}^{2-2\alpha}}\\
&\leq \displaystyle (1+C(\alpha)M)M(t'-t).
\end{array}$$
Then $\theta(t)$ is a Cauchy type at $T^*$. As $ L^2$ is Banach space, then there is an element $\theta^*$ in $L^2$ such that $\theta(t)\rightarrow \theta^*$ in $L^2$ if $t$ goes to $T^*$. As $\theta$ is bounded in the Banach space $H^{2-2\alpha}(\mathbb R^2)$, then $\theta^*\in H^{2-2\alpha}(\mathbb R^2)$.\\
{\bf Step 3:} In this step we prove that $\displaystyle\lim_{t\rightarrow T^*}\|\theta (t)-\theta^*\|_{H^{2-2\alpha}}=0$.\\
By step 2 and by interpolation, since $\theta (t)$ is bounded in $H^{2-2\alpha}$ we have
\begin{equation}\label{eqin}\lim_{t\rightarrow T^*}\|\theta (t)-\theta^*\|_{H^s}=0,\;\forall 0\leq s<2-2\alpha.\end{equation}
Let $0<t< r <T^*$ and a positive real sequence $(s_k)_k$ such  that $1<s_k<2-2\alpha$ and  $s_k\nearrow 2-2\alpha$. We have
$$\|\theta(r)\|^2_{H^{s_k}}+2\int ^{r}_t \| |D|^{\alpha }\theta\|^2_{H^{s_k}}=\|\theta (t)\|^2_{H^{s_k}}-2\underbrace{\int^{r}_t\langle u_{\theta}.\nabla\theta, \theta\rangle_{H^{s_k}}}_{I}.$$
By using Lemma \ref{01}, with $\sigma=s_k$, we get
$$\begin{array}{ll}
|I|&\leq\displaystyle  s_k2^{s_k}C(\alpha)\int_t^{r}\|\theta\|_{\dot {H}^{2-2\alpha}}\|\theta\|_{\dot{H}^{s_k+\alpha}}^2\\
&\leq\displaystyle  (2-2\alpha)2^{2-2\alpha}C(\alpha)\int_t^{r}M\||D|^\alpha\theta\|_{\dot{H}^{s_k}}^2\\
&\leq\displaystyle  (2-2\alpha)2^{2-2\alpha}C(\alpha)M\int_t^{r}\||D|^\alpha\theta\|_{H^{s_k}}^2\\
&\leq\displaystyle  (2-2\alpha)2^{2-2\alpha}C(\alpha)M\int_t^{r}\||D|^\alpha\theta\|_{H^{2-2\alpha}}^2.
\end{array}$$
Then
$$\|\theta(t)\|^2_{H^{s_k}}\leq \|\theta (r)\|^2_{H^{s_k}}+\Big(2(2-2\alpha)2^{2-2\alpha}C(\alpha)M+2\Big)\int ^{r}_t \| |D|^{\alpha }\theta\|^2_{H^{2-2\alpha}}.$$
We pass $r$ to $T^*$, taking into account the equation (\ref{eqin}), we obtain
$$\begin{array}{ll}
\|\theta(t)\|^2_{H^{s_k}}&\leq \displaystyle \|\theta^*\|^2_{H^{s_k}}+\Big(2(2-2\alpha)2^{2-2\alpha}C(\alpha)M+2\Big)\int ^{T^*}_t \| |D|^{\alpha }\theta\|^2_{H^{2-2\alpha}}\\
&\leq \displaystyle \|\theta^*\|^2_{H^{2-2\alpha}}+\Big(2(2-2\alpha)2^{2-2\alpha}C(\alpha)M+2\Big)\int ^{T^*}_t \| |D|^{\alpha }\theta\|^2_{H^{2-2\alpha}}.
\end{array}$$
By Monotone convergence theorem, if $k\rightarrow\infty$, we get
$$\|\theta(t)\|^2_{H^{2-2\alpha}}\leq\|\theta^*\|^2_{H^{2-2\alpha}}+\Big(2(2-2\alpha)2^{2-2\alpha}C(\alpha)M+2\Big)\int ^{T^*}_t \| |D|^{\alpha }\theta\|^2_{H^{2-2\alpha}}.$$
Then, we pass $\displaystyle\limsup_{t\mapsto T^*}$, taking into account the equation (\ref{bleq0}), we obtain $$\limsup_{t\mapsto T^*}\|\theta(t)\|^2_{H^{2-2\alpha}}\leq \|\theta^*\|^2_{H^{2-2\alpha}}.$$
So, by Proposition \ref{thhb} we get $\displaystyle\lim_{t\rightarrow T^*}\|\theta (t)-\theta^*\|_{H^{2-2\alpha}}=0$. Therefore the system of quasi geostrophic starting by $\theta^*$ has a unique solution extends $\theta$, which is absurd.
\section{\bf Global solution in $H^{2-2\alpha}$}
In this section we prove Theorem \ref{TEC2}. Precisely, we prove that if $ \|\theta^0\|_{\dot{H}^{2-2\alpha}}<\frac{1}{4C}$, we get a global solution in $C(\mathbb{R}^+, H^{2-2\alpha}(\mathbb{R}^2))$ satisfying (\ref{eq1th1}). Using Lemma \ref{01} to estimate the non linear part of $(QG)$, we obtain:
$$\|\theta(t)\|^2_{H^{2-2\alpha}}+2 \int^t_0\||D|^{\alpha}\theta\|^2_{H^{2-2\alpha}}\leq \|\theta^0\|^2_{H^{2-2\alpha}}+2C\int^t_0\|\theta\|_{H^{2-2\alpha}}\||D|^{\alpha}\theta\|^2_{H^{2-2\alpha}}.$$
Now, consider the time $T$ defined as follows: $$T= \sup_{0\leq t\leq T^*}\{\sup_{0\leq z\leq t}\|\theta(z)\|_{H^{2-2\alpha}}< 2\|\theta^0\|_{H^{2-2\alpha}} \}.$$
By continuity of the function $(t\rightarrow\|\theta(t)\|_{H^{2-2\alpha}})$ we obtain $T>0$ and for all $0\leq t< T$, we have
\begin{equation}\label{eq2th1}\|\theta(t)\|^2_{H^{2-2\alpha}}+ \int^t_0\||D|^{\alpha}\theta\|^2_{H^{2-2\alpha}}d\tau\leq \|\theta^0\|^2_{H^{2-2\alpha}}.\end{equation}
Again by continuity of the function $(t\rightarrow\|\theta(t)\|_{H^{2-2\alpha}})$ and the above inequality we get $T=T^*$ and $$\int^{T^*}_0\||D|^{\alpha}\theta\|^2_{H^{2-2\alpha}}d\tau<\infty.$$
 Hence $T^*= \infty$, and (\ref{eq1th1}) is given by (\ref{eq2th1}), which complete the proof of Theorem \ref{TEC2}.
\section{\bf Long time decay of global solution}
In this section, we prove Theorem \ref{TEC3}. This proof is done in two steps.\\
\textbf{Step 1}: In this step we shall prove that $\displaystyle\lim_{t\rightarrow\infty} \|\theta(t)\|_{L^2}=0$.\\
We beginning by by recalling the $L^2$ energy estimate:
\begin{equation}\label{eq111}
\|\theta\|_{L^2}^2+2\int_0^t\||D|^\alpha\theta\|_{L^2}^2\leq \|\theta^0\|_{L^2}^2.
\end{equation}
For $\delta>0$, we define the following functions
$$\begin{array}{lcl}
w_\delta= A_\delta(D)\theta&=&\mathcal{F}^{-1}({\bf 1}_{B(0, \delta)}\widehat{\theta}):\;{\rm the\;low\;frequency\;part\;of}\;\;\theta,\\
v_\delta= B_\delta(D)\theta&=& \mathcal{F}^{-1}((1-{\bf 1}_{B(0, \delta)})\widehat{\theta}):\;{\rm the\;high\;frequency\;part\;of}\;\;\theta.
\end{array}$$
and
$$\begin{array}{l}
w_\delta^0= A_\delta(D)\theta^0,\\
v_\delta^0= B_\delta(D)\theta^0.
\end{array}$$
We apply the operator $A_\delta(D)$ to the first equation of system $(QG)$, we obtain
$$ \partial_t w_\delta+|D|^{2\alpha}w_\delta+A_\delta(D)(u_{\theta}.\nabla\theta)=0.$$
Taking the scalar product with $w_\delta$ in last equation and integrating with respect to time, we obtain
$$\| w_\delta(t)\|^2_{L^2}+2\int^t_0\||D|^{\alpha} w_\delta\|^2_{L^2}\leq \|w_\delta ^0\|^2_{L^2} +2\int^t_{0}|\langle A_\delta(D)(u_\theta \nabla \theta), w_\delta\rangle_{L^2}|.$$
Using the fact $A_\delta(D)^2=A_\delta(D)$ and $\langle A_\delta(D)f, g\rangle_{L^2}=\langle f, A_\delta(D)g\rangle_{L^2}$, we get
$$\begin{array}{lcl}
|\langle A_\delta(D)(u_\theta \nabla \theta), w_\delta\rangle_{L^2}|&=&|\langle u_\theta \nabla \theta, w_\delta\rangle_{L^2}|\\
&\leq&\displaystyle\int_{|\xi|<\delta}|\mathcal F({\rm div} (\theta u_\theta))(\xi)|.|\widehat{w}_\delta(\xi)d\xi|\\
&\leq&\displaystyle\int_{|\xi|<\delta}|\xi|.|\mathcal F(\theta u_\theta)(\xi)|.|\widehat{w}_\delta(\xi)d\xi|\\
&\leq&\displaystyle\int_{|\xi|<\delta}|\xi|^{2-2\alpha}.|\xi|^{2\alpha-1}.|\mathcal F(\theta u_\theta)(\xi)|.|\widehat{w}_\delta(\xi)d\xi|\\
&\leq&\displaystyle\delta^{2-2\alpha}\int_{|\xi|<\delta}|\xi|^{2\alpha-1}.|\mathcal F(\theta u_\theta)(\xi)|.|\widehat{w}_\delta(\xi)d\xi|\\
&\leq&\displaystyle\delta^{2-2\alpha}\Big(\int|\xi|^{2(2\alpha-1)}.|\mathcal F(\theta u_\theta)(\xi)|^2d\xi\Big)^{1/2}\|\widehat{w}_\delta\|_{L^2}\\
&\leq&\displaystyle\delta^{2-2\alpha}\|\theta u_\theta\|_{\dot H^{2\alpha-1}}\|\widehat{w}_\delta\|_{L^2}.
\end{array}$$
Applying the product law in homogeneous Sobolev with $s_1=s_2=\alpha$ and using the fact $|\widehat{u_\theta}(\xi)|=|\widehat{\theta}(\xi)|$, we get
$$\|\theta u_\theta\|_{\dot H^{2\alpha-1}}\leq C(\alpha)\|\theta\|_{\dot H^{\alpha}}^2$$
and
$$\begin{array}{lcl}
|\langle A_\delta(D)(u_\theta \nabla \theta), w_\delta\rangle_{L^2}|
&\leq&\displaystyle C(\alpha)\delta^{2-2\alpha}\|\theta\|_{\dot H^\alpha}^2\|w_\delta\|_{L^2}.
\end{array}$$
Using $\|A_\delta(D)f\|_{L^2}\leq \|f\|_{L^2}$ and the $L^2$ energy estimate (\ref{eq111}), we get $\|w_\delta\|_{L^2}\leq \|\theta^0\|_{L^2}$ and
$$\begin{array}{lcl}
\displaystyle\| w_\delta(t)\|^2_{L^2}+2\int^t_0\||D|^{\alpha} w_\delta\|^2_{L^2}&\leq&\displaystyle\|w_\delta ^0\|^2_{L^2}+2C(\alpha)\delta^{2-2\alpha}\|\theta^0\|_{L^2}\int_0^t\|\theta\|_{\dot H^{\alpha}}^2\\
&\leq&\displaystyle\|w_\delta ^0\|^2_{L^2}+C(\alpha)\delta^{2-2\alpha}\|\theta^0\|_{L^2}^3:=\varepsilon_\delta.
\end{array}$$
Clearly $\displaystyle\lim_{\delta\rightarrow0^+}\varepsilon_\delta=0$, then
$$
\lim_{\delta\rightarrow0^+}\sup_{t\geq 0}\| w_\delta\|_{L^2}=0
\;\;{\rm and}\;\;
\lim_{\delta\rightarrow0^+}\int^{\infty}_0 \|w_\delta\|^2_{\dot{H}^{\alpha}}=0.
$$
Let $\varepsilon>0$(suppose that $\varepsilon<1/8C$), then there is $\delta_0>0$ such that: For all $0<\delta\leq \delta_0$, we have
\begin{equation}\label{eqdlta1}
\sup_{t\geq 0}\| w_\delta\|_{L^2}<\varepsilon/2\;\;{\rm and}\;\;\int^{\infty}_0 \|w_\delta\|^2_{\dot{H}^{\alpha}}<\varepsilon/2.
\end{equation}
{\bf In the following we fix a real $\delta\in(0,\delta_0]$}.\\\\
On the other hand
$v_{\delta}$ verifies $$ \partial_t v_{\delta}+|D|^{2\alpha} v_{\delta}+B_{\delta}(D)(u_{\theta}.\nabla\theta)=0,$$then, by Duhamel formula we obtain
$$v_{\delta}=e^{-t|D|^{2\alpha}} v_{\delta}^0-\int_0^te^{-(t-z)|D|^{2\alpha}}B_{\delta}(D)(u_{\theta}.\nabla\theta)dz.$$
Taking the norm in $\dot H^{-\sigma}(\mathbb R^2)$, with $\sigma=2-3\alpha>0$, we get
$$\begin{array}{lcl}
\|v_{\delta}(t)\|_{\dot H^{-\sigma}}&\leq&\displaystyle\|e^{-t|D|^{2\alpha}} v_{\delta}^0\|_{\dot H^{-\sigma}}+\|\int_0^te^{-(t-z)|D|^{2\alpha}}B_{\delta}(D)(u_{\theta}.\nabla\theta)dz\|_{\dot H^{-\sigma}}\\
&\leq&\displaystyle\delta^{-\sigma}e^{-t\delta^{2\alpha}}\|\widehat{v_{\delta}^0}\|_{L^2}+
\Big(\int_{|\xi|>\delta}|\xi|^{-2\sigma}(\int_0^te^{-(t-z)|\xi|^{2\alpha}}|\mathcal F(({\rm div}(\theta u_{\theta})))(z,\xi)|dz)^2d\xi\Big)^{1/2}\\
&\leq&\displaystyle\delta^{-\sigma}e^{-t\delta^{2\alpha}}\|\widehat{\theta^0}\|_{L^2}+
\Big(\int_{|\xi|>\delta}|\xi|^{2-2\sigma}(\int_0^te^{-(t-z)|\xi|^{2\alpha}}|\mathcal F((\theta u_\theta))(z,\xi)|dz)^2d\xi\Big)^{1/2}.
\end{array}$$
Using Lemma \ref{lemf}, we get
$$\|v_{\delta}(t)\|_{\dot H^{-\sigma}}\leq \delta^{-\sigma}e^{-t\delta^{2\alpha}}\|\widehat{\theta^0}\|_{L^2}+\Big(\int_{|\xi|>\delta}|\xi|^{2-2\sigma-2\alpha}\int_0^te^{-(t-z)|\xi|^{2\alpha}}|\mathcal F((\theta u_\theta))(z,\xi)|^2dzd\xi\Big)^{1/2}.$$
Using the fact $2-2\sigma-2\alpha=2(2\alpha-1)$ and the product law in homogeneous Sobolev space with $s_1=s_2=\alpha$, we get
$$\begin{array}{lcl}
\|v_{\delta}(t)\|_{\dot H^{-\sigma}}&\leq& \displaystyle \delta^{-\sigma}e^{-t\delta^{2\alpha}}\|\widehat{\theta^0}\|_{L^2}+\Big(\int_0^te^{-(t-z)\delta^{2\alpha}}\int_{|\xi|>\delta}|\xi|^{2(2\alpha-1)}|\mathcal F((\theta u_\theta))(z,\xi)|^2d\xi dz\Big)^{1/2}\\
&\leq& \displaystyle\delta^{-\sigma}e^{-t\delta^{2\alpha}}\|\widehat{\theta^0}\|_{L^2}+\Big(\int_0^te^{-(t-z)\delta^{2\alpha}}\|(\theta u_\theta)(z)\|_{\dot H^{2\alpha-1}}^2dz\Big)^{1/2}\\
&\leq& \displaystyle\delta^{-\sigma}e^{-t\delta^{2\alpha}}\|\widehat{\theta^0}\|_{L^2}+C(\alpha)\Big(\int_0^te^{-(t-z)\delta^{2\alpha}}\|\theta(z)\|_{\dot H^{\alpha}}^2dz\Big)^{1/2}\\
\end{array}$$
Using the fact
$$\int_0^\infty \Big(\delta^{-\sigma}e^{-t\delta^{2\alpha}}\|\widehat{\theta^0}\|_{L^2}\Big)^2dt=\delta^{-2\sigma}\delta^{-2\alpha}\frac{\|\widehat{\theta^0}\|_{L^2}^2}{2}$$
and
$$\int_0^\infty\Big[\Big(\int_0^te^{-(t-z)\delta^{2\alpha}}\|\theta(z)\|_{\dot H^{\alpha}}^2dz\Big)^{1/2}\Big]^2dt\leq \Big(\int_0^te^{-t\delta^{2\alpha}}dt\Big)\Big(\int_0^\infty\|\theta(t)\|_{\dot H^{\alpha}}^2dt\Big)\leq \delta^{-2\alpha}\|\widehat{\theta^0}\|_{L^2}^2$$
we obtain $v_{\delta}\in L^2(\mathbb R^+,\dot H^{-\sigma}(\mathbb R^2))$ and
$$\int_0^\infty\|v_\delta(t)\|_{\dot H^{-\sigma}}^2dt\leq\Big(\delta^{-2\sigma}\delta^{-2\alpha}\frac{\|\widehat{\theta^0}\|_{L^2}^2}{2}+\delta^{-2\alpha}\|\widehat{\theta^0}\|_{L^2}^2\Big)^2:=M_\delta.$$
The $L^2$ energy estimate (\ref{eq111}) and inequality (\ref{eqdlta1}) give $v_\delta=\theta-w_\delta\in L^2(\mathbb R^+,\dot H^\alpha(\mathbb R^2))$. Therefore $v_{\delta}\in L^2(\mathbb R^+,\dot H^{-\sigma}(\mathbb R^2)\cap \dot H^\alpha(\mathbb R^2))$ and by the injection $$\dot H^{-\sigma}(\mathbb R^2)\cap \dot H^\alpha(\mathbb R^2)\hookrightarrow\dot H^0(\mathbb R^2)=L^2(\mathbb R^2),$$ we get $v_{\delta}\in L^2(\mathbb R^+,L^2(\mathbb R^2))$. Now, we define the set of times
$$E_{\delta}= \{t\geq 0, \|v_{\delta}(t)\|_{L^2}>\frac{\varepsilon}{2}\}.$$ Then
$$ (\frac{\varepsilon}{2})^2 \lambda_1(E_{\delta})\leq\int_{E_{\delta}}\|v_{\delta}(t)\|^2_{L^2}dt\leq \int^\infty _0\|v_{\delta}(t)\|^2_{L^2}dt<\infty.$$
Then, we can define the finite time $T_{\varepsilon}=\displaystyle(\frac{2}{\varepsilon})^2\int^\infty _0\|v_{\delta_0}(t)\|^2_{L^2}dt$ and
$\lambda_1(E_{\delta_0})\leq T_{\varepsilon}$. For $r>0$, there exists $ t_0 \in [0,  T_{\varepsilon}+r] $ such that $t_0 \notin E_{\delta}$ and it results that
\begin{equation}\label{44}
\|v_{\delta}(t_0)\|_{L^2}\leq \frac{\varepsilon}{2}.
\end{equation}
The fact $\theta=v_\delta+w_\delta$ and equations (\ref{eqdlta1}), (\ref{44}) give that
$$\|\theta(t_0)\|_{L^2}\leq \varepsilon.$$
Using the uniqueness of solution ($\gamma(t)=\theta(t_0+t)$) in $H^{2-2\alpha}$ for the following system
$$(QG_{t_0})\hspace {2cm}
   \left\{\begin{array}{l}
   \;\partial_t\gamma+|D|^{2\alpha}\gamma +u_{\gamma}.\nabla \gamma =0\\
   \;\gamma(0, x) =\theta(t_0,x).
   \end{array}\right.$$
we get $$\|\theta(t_0+t)\|_{L^2}^2+2\int_0^t\||D|^\alpha\theta(t_0+z)\|_{L^2}^2dz\leq \|\theta(t_0)\|_{L^2}^2<\varepsilon,\;\;\forall t\geq0.$$
Then $\|\theta(t)\|_{L^2}\rightarrow0$ as $t\rightarrow\infty$, which complete the first step.\\
\textbf{Step 2:} We shall prove that $$\lim_{t\rightarrow}\|\theta(t)\|_{\dot{H}^{2-2\alpha}}=0.$$
As  $\|\theta^0\|_{H^{2-2\alpha}}<\frac{1}{4C}$, then Theorem \ref{TEC2} gives $$\theta\in L^\infty(\mathbb R^+,H^{2-2\alpha}(\mathbb R^2))\cap L^2(\mathbb R^+,\dot H^{2-\alpha}(\mathbb R^2)).$$
Using the interpolation inequality
$$\|\theta(t)\|_{\dot{H}^{2-2\alpha}}\leq \|\theta(t)\|^{\frac{\alpha}{2-\alpha}}_{L^2}\|\theta(t)\|^{\frac{2-2\alpha}{2-\alpha}}_{\dot{H}^{2-\alpha}},$$
we get
$$\theta\in L^{\frac{2-\alpha}{1-\alpha}}(\mathbb R^+,\dot H^{2-2\alpha}(\mathbb R^2)).$$
Let $\varepsilon$ be fixed positif real number and define the set of times
$$F_{\varepsilon}= \{t\geq t_0, \|\theta(t)\|_{\dot H^{2-2\alpha}}\geq\varepsilon\},$$
where $t_0$ is the time given by the first step. Then
$$ \varepsilon^{\frac{2-\alpha}{1-\alpha}} \lambda_1(F_{\varepsilon})\leq\int_{F_{\varepsilon}}\|\theta(t)\|^{^{\frac{2-\alpha}{1-\alpha}}}_{\dot H^{2-2\alpha}}dt\leq \int^\infty _{t_0}\|\theta(t)\|^{^{\frac{2-\alpha}{1-\alpha}}}_{\dot H^{2-2\alpha}}dt<\infty.$$
Then, we can define the finite time $$t_{\varepsilon}=\displaystyle \varepsilon^{-\frac{2-\alpha}{1-\alpha}} \int^\infty _0\|\theta(t)\|^{^{\frac{2-\alpha}{1-\alpha}}}_{\dot H^{2-2\alpha}}dt.$$ We have
$\lambda_1(F_{\varepsilon})\leq t_{\varepsilon}$. Then, there is a time $ t_1 \in [t_0,t_0+t_{\varepsilon}+1] $ such that $t_1 \notin F_{\varepsilon}$ and it results that $\|\theta(t_1)\|_{\dot H^{2-2\alpha}}\leq\varepsilon.$ As $t_1\geq t_0$, then $\|\theta(t_1)\|_{L^2}\leq\varepsilon$ and $\|\theta(t_1)\|_{H^{2-2\alpha}}\leq2\varepsilon.$\\
By choosing $\varepsilon<1/8C$, we get
$$\|\theta(t)\|_{H^{2-2\alpha}}^2+\int_{t_1}^t\||D|^\alpha\theta(z)\|_{H^{2-2\alpha}}^2dz\leq\varepsilon^2,\;\;\forall t\geq t_1.$$
which complete the second step, and the proof of Theorem \ref{TEC3} is finished.\\

\end{document}